\documentclass{myifcolog}

\usepackage{rotating}
\usepackage{array}
\usepackage{booktabs}
\usepackage{enumitem}
\setlist[description]{font=\mdseries\itshape,noitemsep}

\usepackage{amsthm}

\newtheoremstyle{break}%
  {}{}%
  {}{}%
  {\bfseries}{}
  {\newline}{}
  \theoremstyle{break}
  \newtheorem{scheme}{Argumentation Scheme}
\newcommand{\ssch}{\leavevmode \vspace{-\baselineskip}}

\title{Dialogue Types, Argumentation Schemes, and Mathematical Practice:\\
Douglas Walton and Mathematics}
\titlerunning{Walton \& Mathematics}
\addauthor[aberdein@fit.edu]{Andrew Aberdein}{School of Arts \& Communication, Florida Institute of Technology, Melbourne FL}
\authorrunning{Aberdein}
\titlethanks{}

\begin{document}
\maketitle

\begin{abstract}
Douglas Walton's multitudinous contributions to the study of argumentation 
seldom, if ever, directly engage with argumentation in mathematics. Nonetheless, several of the innovations with which he is most closely associated lend themselves to improving our understanding of mathematical arguments. 
I concentrate on two such innovations: dialogue types (\S\ref{sec:types}) and argumentation schemes (\S\ref{sec:schemes}).
I argue that both devices are much more applicable to mathematical reasoning than may be commonly supposed.
\end{abstract}

\section{Dialogue Types}\label{sec:types}
Several decades ago, Douglas Walton proposed a classification of dialogue types: different contexts in which argumentation may arise \cite{Walton89,Walton90}. 
His elegant presentation of the key differences between the most central types (from joint work with Erik C. W. Krabbe) 
is summarized in Table~\ref{types}. 
Dialogue types are distinguished by two main factors: the \emph{initial situation} or circumstances in which the interlocutors find themselves and their \emph{main goal} in pursuing a dialogue. 
Some situations admit more goals than others: if the situation is strongly adversarial, the disputants may be seeking a full determination of the matter at hand, requiring one to \emph{persuade} the other;
or they may need to decide on a course of action and \emph{negotiate} a practical consensus; or they may have little intent beyond airing their respective positions, however quarrelsome or \emph{eristic} the exchange. 
Whereas, if the interlocutors are addressing an open problem where neither has any prior commitments, the last of these goals would be incoherent, but the disputants may still \emph{inquire} into the problem or \emph{deliberate} on how best to act.
And if the interaction arises simply because one party has knowledge the other lacks, only the first sort of outcome makes any sense: an \emph{information-seeking} dialogue.
Thus we arrive at 
six principal dialogue types. However, some of these types may be further subdivided or combined \cite[p.~31]{Walton98}, and the classification is not intended to be exhaustive.

\begin{table}[t]
\begin{center}\small
\caption{Walton and Krabbe's systematic survey of dialogue types \protect\cite[p.~80
]{Walton95}}
\label{types}
\begin{tabular}{m{1.5in}@{\hspace{2em}}*{3}{m{1in}}}
\toprule
{\bf Main Goal}&\multicolumn{3}{c}{\bf Initial Situation}\\
\midrule
&Conflict&Open problem&Unsatisfactory \mbox{spread of} \mbox{information}\\
\midrule
\raggedright Stable agreement/ resolution&{\em Persuasion}&{\em Inquiry}&{\em Information Seeking}\\
\addlinespace
\raggedright Practical settlement/ decision (not) to act&{\em Negotiation}&{\em Deliberation}&N/A\\
\addlinespace
\raggedright Reaching a (provisional) accommodation&{\em Eristic}&N/A&N/A\\
\bottomrule
\end{tabular}
\end{center}
\end{table}

Walton itemizes the goals of the interlocutors, individual or collective, and the potential benefits that may accrue from dialogues of each of the main types in Table~\ref{waltontypes} (taken from \cite[p.~605]{Walton97a}; see also \cite[p.~413
]{Walton90}, \cite[p.~66
]{Walton95}).
Different patterns of argument may be appropriate in different dialogue types: what is reasonable in a negotiation would be improper in a persuasion dialogue; almost anything goes in a quarrel but well-conducted inquiries require respect for procedure, and so forth.
Another important feature of the picture that Walton and Krabbe present is the \emph{dialectical shift}: in the course of a dialogue its type may change \cite[pp.~100 ff.]{Walton95}. This can be a positive development---as a conversation unfolds, its participants can productively shift their attention to different ends. But dialectical shifts can also be troublesome, especially if they go unnoticed by one or more of the participants, leading to the use of argumentative tactics that are now contextually inappropriate.

\begin{table}[t]
\caption{Walton's types of dialogue \protect\cite[p.~605]{Walton97a}}
\label{waltontypes}
\begin{center}\small
\begin{tabular}{*{5}{m{.165\textwidth}}}
\toprule
\textbf{\mbox{Type~of} \mbox{Dialogue}}&\textbf{Initial \mbox{Situation}}&\textbf{Individual \mbox{Goals~of} \mbox{Participants}}&\textbf{Collective \mbox{Goal~of} \mbox{Dialogue}}&\textbf{Benefits}\\
\midrule
Persuasion&\raggedright Difference of opinion&Persuade other party&Resolve difference of opinion&Understand positions\\
\addlinespace
Inquiry&Ignorance&Contribute findings&\raggedright Prove or disprove conjecture&Obtain \mbox{knowledge}\\
\addlinespace
Deliberation&\raggedright Contemplation of future consequences&\raggedright Promote personal goals&\raggedright Act on a thoughtful basis&Formulate personal priorities\\
\addlinespace
Negotiation&\raggedright Conflict of interest&\raggedright Maximize gains (self-interest)&\raggedright Settlement (without undue inequity)&Harmony\\
\addlinespace
Information-Seeking&\raggedright One party lacks information&\raggedright Obtain information&\raggedright Transfer of knowledge&Help in goal activity\\
\addlinespace
\raggedright Quarrel (Eristic)&\raggedright Personal conflict&\raggedright Verbally hit out at and humiliate opponent&Reveal deeper conflict&Vent emotions\\
\addlinespace
Debate&Adversarial&Persuade third party&\raggedright Air strongest arguments for both sides&Spread \mbox{information}\\
\addlinespace
Pedagogical&\raggedright Ignorance of one party&\raggedright Teaching and learning&\raggedright Transfer of knowledge&Reserve \mbox{transfer}\\
\bottomrule
\end{tabular}
\end{center}
\end{table}

Walton does not discuss mathematical dialogues but, in other work, his collaborator Krabbe observes that proofs may occur in several different contexts:
\begin{quote}
\begin{enumerate}[noitemsep]
\item thinking up a proof to convince oneself of the truth of some theorem; 
\item thinking up a proof in dialogue with other people (inquiry dialogue\dots); 
\item presenting a proof to one's fellow discussants in an inquiry dialogue (persuasion dialogue embedded in inquiry dialogue\dots); 
\item presenting a proof to other mathematicians, e.g. by publishing it in a journal (persuasion dialogue\dots); 
\item presenting a proof when teaching (information-seeking and persuasion dialogue) 
\cite[p.~457]{Krabbe08}.
\end{enumerate}
\end{quote}
The primary, 
if not exclusive, concern of 
Krabbe's account (and of my own earlier application of dialogue types to mathematics \cite{Aberdein07}) is with proof.
This may reflect what has been called in another context ``proof chauvinism''---a tendency in philosophers of mathematics to privilege proof over other aspects of mathematical practice \cite{DAlessandro18}. Nonetheless, proof is the aspect of mathematical practice where the applicability of {informal} logic is most unexpected. Hence I shall again begin with proofs.

Are proofs always dialogues or can they be monologues?
The conception of the mathematician as isolated genius has a firm grip on the popular imagination \cite{Hersh10}. 
It is true that mathematicians coauthor papers less than most other scientists, and there are some celebrated examples of solitary endeavour, such as Srinavasa Ramanujan labouring in obscurity or Andrew Wiles's years of solo work prior to his surprise announcement of a proof of Fermat's Last Theorem. Nonetheless, this impression is incomplete at best: Ramanujan only began to fulfil his potential after travelling to Cambridge to collaborate with Hardy and Littlewood \cite{Kanigel91}; Wiles discovered gaps in his solo work which were eventually bridged by a collaboration with Richard Taylor \cite{Mozzochi00}.
On the other hand, as Paul Ernest suggests, 
there are many ways in which mathematics is underpinned by ``symbolically mediated exchanges between persons''---conversations or dialogues:
\begin{quote}
First, the ancient origins as well as various modern systems of proof use dialectical or dialogical reasoning, involving the persuasion of others [see also \cite{DutilhNovaes21}.] \dots\ 
Second, mathematics is a symbolic activity using written inscriptions and language; it inevitably addresses a reader, real or imagined, so mathematical knowledge representations are conversational. 
Third, many mathematical concepts [such as epsilon-delta definitions of limit in analysis and hypothesis testing in statistics] have an internal conversational structure. 
Fourth, the epistemological foundations of mathematical knowledge, including the nature and mechanisms of mathematical knowledge genesis and warranting, utilise the deployment of conversation in an explicitly and constitutively dialectical way. 
Fifth, \dots\ 
mathematical facts stand on the basis of collective agreement and are part of institutional reality \dots\ 
built on interpersonal communicative interactions, that is, through conversation \cite[p.~74]{Ernest16}. 
\end{quote}

Once we agree that mathematical proof is dialogical, we may ask in what dialogue type it characteristically arises.
As Krabbe indicates, the proving process, at least at its inception, might best be thought of as an inquiry dialogue: a collaborative exchange between mathematicians with the shared goal of settling an open question, which neither of them has prejudged.
Certainly such exchanges can be found in mathematics, at least in the context of discovery of mathematical results (see, for example, \cite{Villani15}).
However, there is also an 
unavoidable element of adversariality in the epistemology of mathematical proof:
mathematicians only trust proofs that have gained wide assent from the mathematical audience \cite{Ashton20}; 
the value of that assent lies in the assumption that the proofs have been sufficiently challenged.
Catarina Dutilh Novaes has sought to capture this idea in terms of
prover/sceptic dialogues \cite{DutilhNovaes16,DutilhNovaes18,DutilhNovaes21}.
Prover and sceptic are (idealizations of) two complementary roles in the process that leads to the eventual acceptance of a proof by the mathematical public: the prover presents a putative proof; the sceptic responds with searching but fair questions; through their successive exchanges the proof is improved where necessary and eventually comes to be generally accepted (or is exposed as unsound).
Prover/sceptic dialogues are persuasion dialogues because the parties start from a difference of opinion, as required by their contrasting roles. 

Journal referees can play the sceptic role, 
at least if they are sufficiently thorough in their scrutiny \cite{Andersen20}.
But so can collaborators, at an earlier stage in the development of a proof, or other mathematicians, at a later stage, who expand or refine the published proof in their own work.
Imre Lakatos's celebrated imaginative reconstruction of the development of a proof of the Descartes--Euler conjecture (linking the numbers of vertices, edges, and faces of polyhedra) takes the form of a dialogue between characters loosely representing various nineteenth-century mathematicians \cite{Lakatos76}.
Lakatos identifies a range of dialectical manoeuvres whereby mathematicians either present apparent counterexamples of various kinds to a working conjecture or respond to such apparent counterexamples.
(Alison Pease and colleagues have shown how these Lakatosian manoeuvres can be captured in terms of dialogue games \cite{Pease17}.)


But perhaps the dialogue type 
in which proofs are most frequently presented is neither inquiry nor persuasion, but pedagogical information-seeking. Proofs are presented in countless classrooms at school and university level and even research mathematicians will attend essentially didactic presentations of novel but settled results. Such exchanges are best understood as information-seeking dialogues.
Hence, as Krabbe notes, the development of a proof may be seen as a sequence of dialectical shifts, from an initial inquiry phase, to a more verification-focussed persuasion dialogue, and eventually, if the proof survives these earlier stages, to a dissemination phase characterized by information-seeking dialogues.
Of course, the progress of most significant proofs is seldom this smooth, so the dialectical shifts are likely to be more numerous, as failed attempts at verification send mathematicians back to more open-ended inquiry, or at least open up subsidiary discussions of how localized problems may be addressed.
Michael Barany and Donald MacKenzie, in an ethnographic treatment of mathematical research, describe how some of these processes can work:
\begin{quote}
When a suitable partial result is obtained and researchers are confident in the theoretical soundness of their work, they transition to ``writing up''. Only then do most of the formalisms associated with official mathematics emerge, often with frustrating difficulty. Every researcher interviewed had stories about conclusions that either had come apart in the attempt to formalize them or had been found in error even after the paper had been drafted, submitted, or accepted. Most saw writing-up as a process of verification as much as of presentation, even though they viewed the 
mathematical effort of writing-up as predominantly ``technical'', and thus implicitly not an obstacle to the result's ultimate correctness or insightfulness \cite[pp.~111 f.]{Barany11}.
\end{quote}

Although inquiry, persuasion, and information-seeking dialogues are perhaps the dialogue types most hospitable to proof, they do not exhaust the range of dialogue types in which mathematical argumentation may occur.
By analogy with the device used to indicate problematic sporting records, I have elsewhere used ``proof*'' to refer to ``species of alleged `proof', where there is either no consensus that the method provides proof, or there is broad consensus that it doesn't, but a vocal minority or an historical precedent which points the other way'' \cite[p.~2]{Aberdein09b}.
Amongst the proofs* I included ``proofs* predating modern standards of rigour, picture proofs*, probabilistic proofs*, computer-assisted proofs*, textbook proofs* which are didactically useful but would not satisfy an expert practitioner, and proofs* from neighbouring disciplines with different standards''. 
Each of these cases can be seen two ways: either as a (perhaps very) disputed form of mathematical proof or as an undisputed form of mathematical reasoning that ought to be characterized as something other than proof. Hence, if our goal is to repudiate proof chauvinism and characterize mathematical reasoning in general
, then we must pay attention to proofs*.

\begin{table}[t]
\caption{Some mathematical dialogue types}\label{mathtypes}
\begin{center}\small
\begin{tabular}{*{5}{m{.16\textwidth}}}
\toprule
\textbf{
\mbox{Dialogue} Type}&\textbf{Initial \mbox{Situation}}&\textbf{Main Goal}&\textbf{\mbox{Goal~of} \mbox{Proponent
}}&\textbf{\mbox{Goal~of} \mbox{Respondent
}}\\
\midrule
\raggedright 
Inquiry&Open-mindedness&\raggedright {Prove or} disprove \mbox{conjecture}& \raggedright Contribute to main goal& Obtain \mbox{knowledge}\\
\addlinespace
\raggedright 
Persuasion& \raggedright Difference of opinion& Resolve difference {of opinion with} rigour& \raggedright Persuade respondent
&Persuade \mbox{proponent
}\\
\addlinespace
\raggedright 
Pedagogical Information-Seeking&\raggedright Respondent 
lacks information&\raggedright Transfer of knowledge&\raggedright Disseminate knowledge of results and methods&Obtain \mbox{knowledge}\\
\addlinespace
\raggedright 
Oracular Information-Seeking
& \raggedright 
Proponent lacks information& \raggedright Transfer of knowledge&Obtain \mbox{information}&{Inscrutable}\\
\addlinespace
\raggedright 
Deliberation
&Open-mindedness& \raggedright Reach a provisional conclusion& \raggedright Contribute to main goal& Obtain \mbox{warranted} belief\\
\addlinespace
\raggedright 
Negotiation
& \raggedright Difference of opinion& \raggedright Exchange resources for a provisional conclusion&\raggedright Contribute to main goal&Maximize \mbox{value of} \mbox{exchange}\\
\bottomrule
\end{tabular}
\par\medskip
\end{center}
\end{table}

Table~\ref{mathtypes}, adapted from \cite[p.~148]{Aberdein07}, summarizes the principal mathematical dialogue types discussed so far: inquiry, persuasion, and pedagogical information-seeking.
It also lists three 
dialogue types 
in which proof* is likely to be more at home than proof: deliberation, negotiation, and a non-pedagogical form of information-seeking.
Deliberation and negotiation abandon the goal of stable resolution that we would normally expect of proof whereas oracular information-seeking pursues that goal in an unconventional manner.
In one of his foundational papers on computability, Alan Turing briefly considers the case of a machine ``supplied with some unspecified means of solving number-theoretic problems; a kind of oracle as it were'' \cite[p.~172]{Turing39}.
Subsequent authors expanded this remark into a theory of relative computability \cite{Soare09}. 
There is nothing necessarily supernatural about an oracle machine: a laptop with access to an online database would meet the broad definition (if we ignore Turing's statement that the oracle ``cannot be a machine'' \cite[p.~173]{Turing39}).
However, an oracle is by definition a ``black box'': its inner workings are inscrutable to the local machine; in principle, they could be inscrutable to any analysis.
Sceptics of the proof status of unsurveyably large computer-assisted proofs, such as Thomas Tymoczko, 
have suggested that the appeal such proofs make to a computer should be seen in similar terms \cite{Tymoczko79}.
Analogously, Yehuda Rav proposes as a thought experiment a computer that could answer any mathematical question with certainty but without proof. For Rav, such a machine would be ``a death blow to mathematics, for we would cease having ideas and candidates for conjectures'' \cite[p.~6]{Rav99}. 
Tymoczko and Rav are both concerned about fallacious appeal to authority in mathematical proof, an issue I will return to below.

The 
combinatorialist Edward Swart was also concerned with
``lengthy proofs (whether achieved by hand or on a computer)''. 
He coined the term ``agnograms'' to refer to the resulting ``theoremlike statements'' 
since we are, at least for the immediate future, required to be agnostic about their truth value, as they ``have been neither adequately formalized nor adequately surveyed and are suggestive rather than definitive'', due to the limitations of our available resources \cite[p.~705]{Swart80}.
Establishing an agnogram is thus more of a practical settlement than a stable resolution, suggesting that the dialogue in which it results may better be seen as deliberation or even negotiation, rather than inquiry or persuasion.
Likewise, in a widely discussed polemical proposal, Arthur Jaffe and Frank Quinn sought a clear demarcation between 
``speculative and intuitive work'' in ``theoretical mathematics'' and a ``proof-oriented phase'' of ``rigorous mathematics'' \cite[p.~2]{Jaffe93}.
Of course, speculative and intuitive work is characteristic of the earlier, inquiry phase of a proof dialogue. However, Jaffe and Quinn 
anticipate an outlet for responsibly labelled speculation; since this is provisional in character, the process by which it is derived might be seen as deliberation.
Something similar might also be said about conjectures, particularly the wide-ranging, fruitful conjectures that comprise the framework of mathematical research programmes, sometimes called ``architectural conjectures'' \cite[p.~198]{Mazur97}.
Even more speculative is the suggestion of Doron Zeilberger that in the not so distant future ``semi-rigorous mathematics'' may essentially assign price tickets to proofs, indicating the quantity of computational resources needed for certainty, thereby situating mathematical proof within a negotiation dialogue
\cite{Zeilberger93}. 
This proposal has not generally been well received \cite{Andrews94}. Nonetheless, in a weaker form it reflects a truism: even the purest of mathematicians cannot ignore issues of funding, even if the link to their work is not as intimate as Zeilberger suggests. 
Lastly, even eristic dialogues have had some role to play in mathematical reasoning, as witnessed by such celebrated quarrels as that between the early modern mathematicians Girolamo Cardano and Niccol\`o Tartaglia \cite{Rothman14}.
The salient detail is not the asperity of their exchange, which ultimately turned on an accusation of theft of intellectual property, but the adversarial strategy mathematicians of that era adopted to convince the mathematical public of their successes. Rival mathematicians would keep their methods (in this case of solving cubic equations) secret but challenge each other to public contests, each solving problems set by the other until the winner posed a problem the loser could not solve.


\section{Argumentation Schemes}\label{sec:schemes}
An argumentation scheme is a stereotypical pattern of reasoning. 
In recent decades, the study and classification of argumentation schemes has been the most influential aspect of Douglas Walton's work \cite{Walton96,Walton08}. 
Although antecedents of the argumentation scheme can be traced back millennia to the tradition of loci or topoi, Walton's work set it on a new foundation of rigour and clarity.
Building on that foundation, Hans Hansen has proposed the following definition of argumentation scheme:
``(i) a pattern of argument, (ii) made of a sequence of sentential forms with variables, of which (iii) at least one of the sentential forms contains a use of a schematic constant or a use of a schematic quantifier, and (iv) the last sentential form is introduced by a conclusion indicator like `so' or `therefore'\,'' \cite[p.~349]{Hansen20a}.
Schemes also generally include `critical questions', which itemize possible lines of response. The critical questions are key to the evaluation of defeasible schemes: whether the argument should be judged to have succeeded or whether it has been defeated will turn on whether the questions can receive a satisfactory answer.

Walton argues that in principle all defeasible argumentation schemes could be understood as special cases of a defeasible version of modus ponens \cite[p.~366]{Walton08}:
\begin{scheme}
[Defeasible Modus Ponens]
\label{DMP}\ssch
\begin{description}
\item[Data:]$P$.
\item[Warrant:] As a rule, if $P$, then $Q$.\\
Therefore, \dots 
\item[Qualifier:] presumably, \dots
\item[Claim
:] \dots\ $Q$. 
\end{description}
\begin{center}
Critical Questions:
\end{center}
\begin{enumerate}[noitemsep]
\item \emph{Backing:} What reason is there to accept that, as a rule, if $P$, then $Q$?
\item \emph{Rebuttal:} Is the present case an exception to the rule that if $P$, then $Q$? 
\end{enumerate}
\end{scheme}
I have reconstructed Walton's scheme for defeasible modus ponens so as to bring out its resemblance to another very general model of defeasible reasoning, the Toulmin layout \cite[p.~829]{Aberdein19a}.
(On the relationship of schemes to layouts, see also \cite[pp.~22 ff.]{Pease11}; for a contrasting view, see \cite{Hitchcock03}.)
This is not an accidental choice: Toulmin layouts have lately found widespread employment in the analysis of mathematical argumentation (for recent surveys, see \cite{Knipping19,Krummheuer15}).
Although Walton emphasized defeasible schemes, deductive rules of inference can also be seen as argumentation schemes: 
the schemes framework is ``illatively neutral'' \cite[p.~355]{Hansen20a}.
This means that argumentation schemes can provide a unified treatment of a wide range of arguments employed in mathematics.
Indeed, a number of authors have applied argumentation schemes to mathematical reasoning
\cite{Aberdein07b,Aberdein09,Aberdein13,Aberdein13a,Aberdein19a,Dove09,Metaxas15,Metaxas16,Pease11}.

The illative neutrality 
of the schemes framework
licences scepticism about the ``standard view'' \cite{Marfori10} of mathematical argumentation as purely comprised of derivations, that is arguments in which every step instantiates a deductive inference rule.
To that end, I have elsewhere proposed a threefold distinction between A-, B-, and C-schemes: 
\begin{quote}
\begin{itemize}
\item {\bf A-schemes} {correspond directly to derivation rules}. 
(Equivalently, we could think in terms of a single A-scheme, the `pointing scheme' which picks out a  derivation whose premisses and conclusion are formal counterparts of its data and claim.) 
\item {\bf B-schemes} are {exclusively mathematical arguments}: high-level algorithms or macros. Their instantiations {correspond to substructures of derivations} rather than individual derivations (and they may appeal to additional formally verified propositions). 
\item {\bf C-schemes} are even looser in their relationship to derivations, since {the link between their data and claim need not be deductive}. 
Specific instantiations may still correspond to derivations, but there will be no guarantee that this is so and no procedure that will always yield the required structure even when it exists.
Thus, where the qualifier of A- and B-schemes will always indicate deductive certainty, the qualifiers of C-schemes may exhibit more diversity. Indeed, different instantiations of the same scheme may have different qualifiers (\cite[p.~829]{Aberdein19a}; cf.~\cite[pp.~366 f.]{Aberdein13a}). 
\end{itemize}
\end{quote}
So the widespread 
``standard'' view of mathematical proof, that it is identical to derivation, could be expressed as denying C-schemes a place in proofs. I have argued against that view \cite[p.~375]{Aberdein13a}, but even if it were to be conceded, it would still leave room for C-schemes in other forms of mathematical reasoning.

What sort of schemes might C-schemes be? Some of them may be unique to mathematics, but we should expect others to resemble schemes that have been found useful in addressing non-mathematical reasoning.
Walton and his collaborators have made a number of attempts to classify such general purpose argumentation schemes. 
Table~\ref{schemes} is based on a recent classification he developed with Fabrizio Macagno. 
Walton and Macagno employ a series of binary distinctions: first between source-dependent arguments and source-independent arguments; then subdividing the latter into practical reasoning and epistemic reasoning; which is in turn divided into discovery arguments and arguments applying rules to cases. Each of the four resulting headings are then further subdivided into various thematic groups of individual schemes.
However, Walton and Macagno concede that this classification is incomplete, notably omitting some linguistic arguments \cite[p.~24]{Walton16a}.
  
\begin{sidewaystable}[htp]
\begin{center}\tiny
\begin{tabular}{*{3}{p{.22\textwidth}}p{.25\textwidth}}
\toprule
Discovery arguments&Applying rules to cases&Practical reasoning&Source-dependent arguments\\
\midrule
\begin{enumerate}[leftmargin=1.5em]
\item Arguments establishing rules
\begin{itemize}
\item Argument from a random sample to a population
\item Argument from best explanation
\end{itemize}
\item Arguments finding entities
\begin{itemize}
\item Argument from sign \cite{Dove09,Metaxas15,Metaxas16}
\item Argument from ignorance \cite{Aberdein13}
\end{itemize}
\end{enumerate}\vfill
&
\begin{enumerate}[leftmargin=1.5em]
\item Arguments based on cases
\begin{itemize}
\item Argument from an established rule
\item\raggedright Argument from verbal classification \cite{Aberdein09,Metaxas15,Metaxas16}
\item Argument from cause to effect
\end{itemize}
\item Defeasible rule-based arguments
\begin{itemize}
\item Argument from example \cite{Aberdein13,Aberdein19a,Metaxas15,Metaxas16}
\item Argument from analogy \cite{Aberdein13a,Metaxas16}
\item Argument from precedent \cite{Pease11}
\end{itemize}
\item Chained arguments connecting rules and cases
\begin{itemize}
\item Argument from gradualism \cite{Aberdein13}
\item Precedent slippery slope argument
\item Sorites slippery slope argument
\end{itemize}
\end{enumerate}&
\begin{enumerate}[leftmargin=1.5em]
\item\raggedright Instrumental argument from practical reasoning
\begin{itemize}
\item Argument from action to motive
\end{itemize}
\item Argument from values
\begin{itemize}
\item Argument from fairness
\end{itemize}
\item Value-based argument from practical reasoning
\begin{enumerate}
\item Argument from positive or negative consequences \cite{Aberdein13,Metaxas15,Metaxas16}
\begin{itemize}
\item Argument from waste
\item Argument from threat
\item Argument from sunk costs
\end{itemize}
\end{enumerate}
\end{enumerate}
&
\begin{enumerate}[leftmargin=1.5em]
\item Arguments from position to know 
\begin{enumerate}
\item\raggedright Argument from expert opinion \cite{Aberdein07b,Aberdein13,Metaxas15}
\item Argument from position to know \cite{Aberdein13}
\begin{itemize}
\item Argument from witness testimony
\end{itemize}
\end{enumerate}
\item Ad hominem arguments
\begin{enumerate}
\item Direct ad hominem 
\item\raggedright Circumstantial ad hominem 
\begin{itemize}
\item Argument from inconsistent commitment
\item Arguments attacking personal credibility
\begin{enumerate}
\item Arguments from allegation of bias
\item Poisoning the well by alleging group bias
\end{enumerate}
\end{itemize}
\end{enumerate}
\item Arguments from popular acceptance
\begin{itemize}
\item\raggedright Argument from popular opinion \cite{Aberdein09}
\item Argument from popular practice \cite{Aberdein09}
\end{itemize}
\end{enumerate}
\\
\bottomrule
\end{tabular}
\end{center}
\caption{Walton \& Macagno's 
partial classification of schemes (adapted from \cite[p.~22
]{Walton16a}), with 
applications to mathematical argumentation indicated.}
\label{schemes}
\end{sidewaystable}%

I have annotated Table~\ref{schemes} with citations to works in which mathematical versions of each scheme are discussed. As may be seen, mathematical arguments have been identified under each of Walton and Macagno's four main headings. 
In addition, mathematical applications have been found for several schemes that are missing from Table~\ref{schemes} but which are found in the more exhaustive (but less structured) list in \cite{Walton08}.
These include linguistic arguments, such as arguments from arbitrariness or vagueness of a verbal classification \cite{Pease11} and argument from definition to verbal classification \cite{Aberdein13}, but also  
source-dependent arguments, such as ethotic argument \cite{Aberdein13}, 
practical reasoning arguments, such as argument from 
positive consequences \cite{Aberdein13}, 
discovery arguments, such as abductive argument \cite{Pease11} and argument from evidence to a hypothesis \cite{Aberdein13,Aberdein19a,Pease11}, and 
arguments applying rules to cases, such as argument from an exceptional case \cite{Pease11}.
Conversely, not all of the individual schemes 
in Table~\ref{schemes} have yet been found useful in discussing mathematics.
While some of these omissions may merely be oversights, others are to be expected.
For example, causal reasoning, whether argument from cause to effect or the various kinds of slippery slope, is unlikely to be of direct application to mathematics, since mathematical objects are generally understood to be causally inert.
In the remainder of this section, I will discuss the mathematical applications of a sample of schemes, chosen in part to remedy some of the omissions in Table~\ref{schemes}.

\subsection{Epistemic reasoning}
Walton and Macagno subdivide epistemic reasoning into two subcategories, discovery arguments and arguments applying rules to cases.
Arguments of both kinds can be readily found in mathematical reasoning.
In particular, many discovery arguments are broadly abductive in character and abduction has been proposed as an account of mathematical reasoning in a wide range of situations, including 
classroom discussion \cite{Ferrando05,Meyer10};
concept formation in mathematical practice \cite{Heeffer07}; and
the selection and defence of axioms \cite{Heron20a}.
There are several abductive schemes in Walton's catalogue, 
including multiple subtypes of abductive argument \cite[p.~329]{Walton08} and 
argument from evidence to (verification of) a hypothesis, which I have discussed elsewhere \cite{Aberdein13,Aberdein19a,Pease11}. 
Another such scheme is argument from sign:
\begin{scheme}[
Argument from Sign {\cite[p.~329]{Walton08}}]
\ssch\label{sign}
\begin{description}
\item[Specific Premise:] $A$ (a finding) is true in this situation.
\item[General Premise:] $B$ is generally indicated as true when its sign, $A$, is true.
\item[Conclusion:] 
$B$ is true in this situation.
\end{description}
\begin{center}
Critical Questions:
\end{center}
\begin{enumerate}[noitemsep]
\item What is the strength of the correlation of the sign with the event signified?
\item Are there other events that would more reliably account for the sign?
\end{enumerate}
\end{scheme}
Argument from sign has been discussed since antiquity, particularly in the context of medical reasoning \cite{Allen01}.
The explicit application of Scheme~\ref{sign} to mathematics is due to Ian Dove, who uses it to analyse a surprising but widely discussed class of proofs*: those employing molecular computation \cite{Dove09}.
This consists in encoding 
a mathematical problem 
into strands of DNA which are then subject to standard laboratory assays that determine the solution of the problem with high likelihood \cite{Adleman94}. Hence the outcome of the assay is a sign of the mathematical problem having a specific solution, and the mathematician infers the latter from the former in accordance with Scheme~\ref{sign}.
This, and other less esoteric probabilistic methods, such as the Miller--Rabin primality test, are generally viewed by mathematicians as heuristically useful but falling short of the standards of rigour required for proof. Nonetheless, the intellectual defensibility of this perspective has also been the subject of a debate in philosophy of mathematics \cite{Easwaran09,Fallis97,Fallis11}.
Much of this debate could be understood as offering competing answers to the critical questions for Scheme~\ref{sign}.
Dove also suggests that what I have referred to above as oracular information seeking could be analysed as employing the same scheme \cite[p.~144]{Dove09}.

Argument from sign also illustrates the importance of the illative neutrality of argumentation schemes: not all its instances need be defeasible. We can find deductive instances of Scheme~\ref{sign}.
For example, much of twentieth and twenty-first century mathematics employs increasingly complex mathematical infrastructure or tools: that is, mathematical theories designed to help us investigate other areas of mathematics. Mathematical tools, such as Galois theory or \emph{K}-theory, establish rigorous relationships between outwardly unrelated classes of mathematical objects. As Jean-Pierre Marquis observes, the function of such tools is to
``reveal important \emph{properties} of the objects studied, and \emph{only} these properties'' \cite[p.~264]{Marquis97}.
In other words, a result in one of the two related areas may be taken as a sign that one of the presumably less tractable objects in the other area has a particular property. However, since the relationship between the areas can be rigorously established, the sign is not merely generally indicative, but infallibly so.

Applying rules to cases is also a very widespread practice in mathematics. Several of the schemes that fall under this heading, such as argument from verbal classification, argument from example, and argument from analogy, have mathematical applications that I have discussed elsewhere \cite{Aberdein09,Aberdein13,Aberdein13a}. Walton and Macagno also include chained arguments, which comprise a substantial proportion of mathematical reasoning \cite[p.~235]{Aberdein13}.
But here I shall focus on a different scheme:
\begin{scheme}[
Argument from an Established Rule \protect{\cite[p.~343]{Walton08}}]
\ssch\label{rule}
\begin{description}
\item[Major Premise:] If  
carrying out types of actions including $A$ is the established rule for $x$, then (unless the case is an exception), $x$ must carry out $A$.
\item[Minor Premise:] Carrying out types of actions including $A$ is the established rule for $a$.
\item[Conclusion:] Therefore, $a$ must carry out $A$.
\end{description}
\begin{center}
Critical Questions:
\end{center}
\begin{enumerate}[noitemsep]
\item Does the rule require carrying out types of actions that include $A$ as an instance?
\item Are there other established rules that might conflict with or override this one?
\item Is this case an exceptional one, that is, could there be extenuating circumstances or an excuse for noncompliance?
\end{enumerate}
\end{scheme}
Scheme~\ref{rule} is framed in terms of actions to be carried out. That might initially appear to be an obstacle to its application to mathematics. 
However, as Wilfrid Hodges has observed, 
informal mathematical arguments 
include not only
the ``object sentences'', in which some mathematical content is explicitly given, and 
``stated or implied justifications for putting the object sentences in the places where they appear'' 
but also ``instructions to do certain things which are needed for the proof'' \cite[p.~6]{Hodges98}.
The last of these, carrying out actions, has perhaps received least attention from logicians, but it is ubiquitous and important. Many proofs instruct us to ``\,`Suppose $C$', `Draw the following picture, and consider the circles $D$ and $E$', `Define $F$ as follows'\,'' and so forth \cite[p.~6]{Hodges98}. 
Language of this sort is phrased conventionally as an instruction to the reader, but it is also a description of the actions undertaken by the deviser of the proof. But how did the proof's author know which actions to carry out? At least in some cases, by application of Scheme~\ref{rule}.

Carrying out a rule has also been the focus of a significant debate in the philosophy of mathematics, inspired by the work of Ludwig Wittgenstein. 
Wittgenstein considers the case of a pupil who learns to follow a rule whereby he writes down a series of natural numbers each greater than its predecessor by $2$. However, after he gets to $1000$ he increments the numbers by $4$, instead of $2$, but takes himself still to be following the same rule \cite[\S185]{Wittgenstein53}. What ought we to make of such behaviour?
It has been suggested that Wittgenstein's intent was to suggest a general scepticism about rule-following \cite{Kripke82}. If that were to be the case, Critical Question 2 in Scheme~\ref{rule} would always receive an affirmative answer: there would always be another rule which might override any rule we may consider.
Less radically, we could read Wittgenstein as counselling against a platonist interpretation of rules as existing independently of the practices they govern \cite[p.~91]{Wright90}. Rather we should understand rules as implicit within our practice but nonetheless as carrying normative force. The ontological status of rules is the subject of a difficult and important debate.
Fortunately, Scheme~\ref{rule}, and related rule-establishing and applying schemes, are neutral as to the outcome of that debate.

\subsection{Practical reasoning}
Practical reasoning is an inevitable component of 
resource-sensitive mathematical deliberation dialogues 
whether limited by time, money, or processor capacity. 
If numerical approximation methods are easy and cheap 
and an exact answer would be expensive and slow, we may settle for the former.
More broadly, a dialogue can shift to addressing this sort of question 
within reasoning about a problem whenever a choice of methods arises. 
So practical reasoning is not 
just a project management phase to be completed before the real work begins, but potentially a recurrent phenomenon throughout the research process. 
For example, James Franklin points out one context in which practical reasoning occurs in the career of almost every research mathematician: choice of Ph.D.\ topic. A Ph.D.\ thesis is expected to address an open question which must also be
``tractable, that is, probably solvable, or at least partially solvable, by three years' work at the Ph.D.\ level'' \cite[p.~2]{Franklin87}.
Determining whether a problem is tractable is not something which can be established with certainty. But it is critical to the success of the Ph.D.\ 
Elsewhere I have addressed some special cases of practical reasoning, including
argument from positive 
consequences: if important results would follow from a conjecture, that at least provides good reason to investigate it more thoroughly than similar, but less consequential conjectures \cite[pp.~235 f.]{Aberdein13}.
However, I have not directly discussed the most general practical reasoning scheme in Walton's taxonomy 
(see also \cite[p.~131]{Walton95a}):

\begin{scheme}[Practical Inference {\cite[p.~323]{Walton08}}]\ssch\label{prac}
\begin{description}
\item[Major Premise:] I have a goal $G$.
\item[Minor Premise:] Carrying out this action $A$ is a means to realise $G$.
\item[Conclusion:] Therefore, I ought (practically speaking) to carry out this action $A$.
\end{description}
\begin{center}
Critical Questions:
\end{center}
\begin{enumerate}[noitemsep]
\item What other goals that I have that might conflict with $G$ should be considered?
\item What alternative actions to my bringing about $A$ that would also bring about $G$ should be considered?
\item Among bringing about $A$ and these alternative actions, which is arguably the most efficient?
\item What grounds are there for arguing it is practically possible for me to bring about $A$?
\item What consequences of my bringing about $A$ should also be taken into account?
\end{enumerate}
\end{scheme}

Scheme~\ref{prac} could be used to analyse much of the embedded negotiations about resource allocation discussed above. 
However, it can also play a more direct role in mathematical reasoning:
Yacin Hamami and Rebecca Morris have proposed an account of plans and planning in proving in terms of intentions and practical reasoning, 
building on Michael Bratman's work in the philosophy of action \cite{Bratman87}.
The process of finding a proof, at least if the proof is of any complexity, may involve the construction and execution of a carefully devised proof plan, which Hamami and Morris define as
``an ordered tree whose nodes are proving intentions, whose root is the proving intention corresponding to the theorem at hand, and where each set of ordered children consists of a subplan obtained from the parent node through an instance of practical reasoning'' \cite{Hamami20a}.
The plan is not the proof, any more than the map is the journey. 
But, they suggest, the plan is essential not only to successfully finding the proof, but also to subsequently understanding the proof.

\subsection{Source-dependent arguments}
Mathematicians have an ambivalent attitude to authority: there is 
``a schism in the mathematical community \dots\ 
[between those] who think that one should never use a result without having understood its entire proof \dots\ 
[and those who] don’t share that view'' (anonymous mathematician, interviewed in \cite{Andersen21}).
Unlike the empirical sciences, where replication of experiments can require substantial resources, it is in principle always possible for mathematicians to work through every step of every proof they use. 
But, for many mathematicians, a division of labour is unavoidable and even welcome.
Hence there is a place in mathematics for one of the most discussed argumentation schemes, that for argument from expert opinion \cite{Aberdein07b,Aberdein13,Metaxas15}.
Notably, there are several disputes over the role of testimony in mathematical reasoning that may be understood in terms of the critical questions of this scheme.
``Folk theorems'' are one such troublesome case. These are results which are widely used and accepted despite lacking a clear source in the literature.
In a pioneering study drawing attention to their prevalence, the theoretical computer scientist 
David Harel suggests that
``popularity, anonymous authorship, and age \dots\ seem to be necessary and sufficient for a theorem to be folklore, [although] the ways in which they appear and can be established are by no means clear-cut'' \cite[p.~379 f.]{Harel80}.
Don Fallis raises the concern that the citation of folk theorems may represent ``universally untraversed gaps''  in mathematical reasoning since ``everyone is convinced that these theorems are provable, but no one has bothered to work through all the details of a proof'' \cite[p.~62]{Fallis03}. And, of course, everyone could be wrong.
Colin Rittberg and colleagues raise a different problem: the ambiguous status of folk theorems, neither rigorously proved nor strictly open problems, presents a hazard to young researchers \cite{Rittberg19}. Actually proving a folk theorem can be an unrewarding project, since referees may reject such work as unoriginal---despite being unable to cite any prior proof.
A possible resolution to this and related problems lies in the work of
Kenny Easwaran, who has posited a property of ``transferability'' that may distinguish the proofs which safely support argument from expert opinion from those that do not. Transferable proofs are those that
``rely only on premises that the competent reader can be assumed to antecedently believe, and only make inferences that the competent reader would be expected to accept on her own consideration'' \cite[p.~354]{Easwaran09}.
Hence folk theorems lack transferable proofs, unless there is a proof simple enough for any competent mathematician to reconstruct. 
Many other proofs* would be untransferable too, including unsurveyably long proofs, probabilistic proofs, and proofs that rely on empirical procedures.
This suggests a revision, or precisification, of the critical questions of the expert opinion scheme.

Argument from expert opinion is not the only source-dependent argument relevant to mathematics.
Walton draws an important distinction between argument from expert opinion and argument from position to know. 
The distinction is familiar from legal practice, as that between expert and fact witnesses.
I have suggested elsewhere that argument from position to know provides a model for appeals to intuition in mathematics \cite[p.~240 f.]{Aberdein13}.
In this I follow philosophers, such as Elijah Chudnoff, for whom intuition is analogous to perception 
\cite{Chudnoff13c,Chudnoff13}.
Reports on (reliable) perceptions support cogent instances of argument from position to know, so if intuition may be treated analogously, then the same scheme should apply.
Arguments from popular opinion and popular practice also have important applications to mathematics \cite[p.~283 ff.]{Aberdein09}. 
However, here I will focus on yet another source-dependent argument:
\begin{scheme}[Ethotic Argument {\cite[p.~336]{Walton08}}]\label{ethotic}
\leavevmode \vspace{-\baselineskip}
\begin{description}
\item[Major Premise:] If $x$ is a person of good (bad) moral character, then what $x$ says 
should be accepted as more plausible (rejected as less plausible).
\item[Minor Premise:] $a$ is a person of good (bad) moral character.
\item[Conclusion:] Therefore, what $a$ says 
should be accepted as more plausible (rejected as less plausible).
\end{description}
\begin{center}
Critical Questions:
\end{center}
\begin{enumerate}[noitemsep]
\item Is $a$ a person of good (bad) moral character?
\item Is character relevant in the dialogue?
\item Is the weight of presumption claimed strongly enough warranted by the evidence given? 
\end{enumerate}
\end{scheme}
Superficially, this scheme may appear to be a poor fit for mathematical argument, but there is empirical research that suggests the applicability of something much like it. 
Matthew Inglis and Juan Pablo Mej\'ia-Ramos 
gave an informal mathematical argument for the presence of one million sevens in the decimal expansion of $\pi$ to samples of undergraduates and research mathematicians and asked them to rate how persuasive they found it \cite{Inglis09a}. 
Participants in both groups for whom the argument was correctly attributed to 
the prominent mathematician Tim Gowers ranked it as more persuasive than those for whom it was presented anonymously, significantly so for the researchers (who were presumably more likely to have heard of Gowers).
Of course, it is not Gowers's (doubtless exemplary) moral conduct which leads us to trust his arguments, but rather his demonstrably high standards as a working mathematician. As one of the research subjects in this study comments, ``We are told the argument is made by a reputable mathematician, so we implicitly assume  that he would tell us if he knew of any evidence or convincing arguments to the contrary'' \cite[p.~42]{Inglis09a}.
This suggests that we should localize Scheme~\ref{ethotic} to mathematics, 
replacing instances of ``moral'' with ``mathematical'':
\begin{scheme}[Ethotic Mathematical Argument]\label{ethmath}
\leavevmode \vspace{-\baselineskip}
\begin{description}
\item[Major Premise:] If $x$ is a person of good (bad) {mathematical} character, then what $x$ says 
should be accepted as more plausible (rejected as less plausible).
\item[Minor Premise:] $a$ is a person of good (bad) {mathematical} character.
\item[Conclusion:] Therefore, what $a$ says 
should be accepted as more plausible (rejected as less plausible).
\end{description}
\begin{center}
Critical Questions:
\end{center}
\begin{enumerate}[noitemsep]
\item Is $a$ a person of good (bad) {mathematical} character?
\item Is mathematical character relevant in the dialogue?
\item Is the weight of presumption claimed strongly enough warranted by the evidence given? 
\end{enumerate}
\end{scheme}
A (presumably implicit) invocation of Scheme~\ref{ethmath} would explain Inglis and Mej\'ia-Ramos's finding, but it raises other questions, most centrally: what is {mathematical} character?
This is a question which recent work applying virtue epistemology to mathematical practice has sought to answer \cite{Tanswell20}.
Mathematicians have also employed virtue talk to describe their activities. For example,
George P\'olya asserts that the following ``moral qualities'' are required of a mathematician:
\begin{quote}
\begin{itemize}
\item First, we should be {ready to revise} any one of our beliefs.
\item Second, we should {change} a belief when there is a {compelling reason} to change it.
\item Third, we should {not change} a belief {wantonly}, without some good reason \cite[vol.~1, p.~8]{Polya54}.
\end{itemize}
\end{quote}
He goes on to expand on these points in explicitly virtue-theoretic terms, telling us that intellectual courage is required for the first, intellectual honesty for the second, and wise restraint for the last.
The first of these is a widely discussed intellectual virtue and the second is, at least in this context, closely related to the even more widely discussed intellectual humility. Wise restraint is perhaps more familiar as prudence or practical wisdom.
More recent mathematicians who have discussed character virtues relevant to their profession include Michael Harris \cite{Harris15} and Francis Su \cite{Su19}.
Of course, these mathematicians would not necessarily endorse every application of Scheme~\ref{ethmath}. 
Indeed, as we saw above, some mathematicians consider it to be a virtue to never take mathematical results on trust and insist on convincing themselves of the proof of any result they cite.
Nonetheless, many other mathematicians, especially when reasoning speculatively rather than writing up proofs, rely on a division of labour which makes essential use of the informal arguments of their peers, and may be expected to take those arguments more seriously when they have more reason to trust their authors.

\section{Conclusion}
Recent work in the philosophy of mathematical practice has drawn attention to mathematical reasoning in contexts other than proof and
challenged the traditional conception that mathematical proof is essentially reducible to formal derivation. 
This leaves a conspicuous lacuna in our understanding of how mathematics works. Formal logic is an excellent tool for the analysis of formal derivations, but it is less well adapted to the analysis of informal reasoning.
However, the tools developed by informal logicians such as Douglas Walton are a rich source for remedying this deficit.
In particular, as we have seen, dialogue types help to contextualize the different levels of rigour that mathematical argument can exhibit and argumentation schemes provide a valuable taxonomy of the steps that comprise such arguments.

\bibliographystyle{plain}\bibliography{ArgMath}
\end{document}